\documentstyle[12pt,amssymb,amstex]{amsart}
\numberwithin{equation}{section}
\def\ca{{\cal A}}
\def\cb{{\cal B}}

\def\ce{{\cal E}}
\def\cf{{\cal F}}

\def\ch{{\cal H}}

\def\cs{{\cal S}}

\def\ga{{\frak A}}
\def\gb{{\frak B}}

\def\gd{{\frak D}}

\def\bc{{\mathbb C}}

\def\bbf{{\mathbb F}}

\def\bn{{\mathbb N}}

\def\bz{{\mathbb Z}}

\def\a{\alpha}
\def\b{\beta}
  \def\G{\Gamma}
\def\d{\delta}  

\def\e{\epsilon}
\def\l{\lambda} 

\def\m{\mu}

\def\n{\nu}

\def\s{\sigma} 
\def\t{\tau}
\def\f{\varphi} 
\def\th{\theta}  
\def\om{\omega} \def\Om{\Omega}

\newtheorem{thm}{Theorem}[section]
\newtheorem{lem}[thm]{Lemma}
\newtheorem{cor}[thm]{Corollary}
\newtheorem{prop}[thm]{Proposition}
\newtheorem{defin}[thm]{Definition}

\def\di{\mathop{\rm d}}

\def\id{\mathop{\rm id}}
\def\supp{\mathop{\rm supp}}
\def\tr{\mathop{\rm Tr}}
\def\re{\mathop{\rm Re}}

\def\idd{{\bf 1}\!\!{\rm I}}

\begin{document}

\title[ergodic properties]
{on strong ergodic properties of quantum dynamical systems}
\author{Francesco Fidaleo}\footnote{Permanent address:
Dipartimento di Matematica,
II Universit\`{a} di Roma ``Tor Vergata'',
Via della Ricerca Scientifica, 00133 Roma, Italy.
e--mail: \tt fidaleo@@mat.uniroma2.it}
\address{Francesco Fidaleo\\
Department of Mathematics\\
Texas A\&M University \\
Milner Hall, College Station TX 77843--3368, USA}

\begin{abstract}
We show that the the shift on the reduced $C^{*}$--algebras of RD--groups,
including the free
group on infinitely many generators, and the amalgamated free product
$C^{*}$--algebras, enjoys the very strong
ergodic property of the convergence to the equilibrium. Namely, the free shift converges, pointwise in the weak topology, to the conditional expectation onto the fixed--point subalgebra.
Provided the invariant state is unique, we also show that such an ergodic property cannot be
fulfilled by any
classical dynamical system, unless it is conjugate to the trivial one--point dynamical system.
\end{abstract}

\maketitle

\section{introduction}

The study of the ergodic properties of quantum dynamical systems has 
been an impetuos growth in the last years, in view of natural applications to various field of mathematics and physics. It is then of interest to understand among the various 
ergodic properties, which one survives and/or is meaningful, by passing from the classical 
to the quantum case. 
By coming back to the classical case, a very strong ergodic property for the dynamical system 
$(\Om,T)$
consisting of  a compact Hausdorff space $\Om$ and a homeomorphism
$T$, is the
unique ergodicity. This means that there would exist a unique
invariant Borel measure $\m$ for $T$. 
It is seen that  the ergodic
average ${\displaystyle\frac1n\sum_{k=0}^{n-1}f\circ T^{k}}$ of any fixed function $f$,
converges uniformly to the constant function $\int f\di\m$. A pivotal example of a classical uniquely ergodic dynamical
system is given by an irrational rotation on the unit circle. In the quantum setting, the unique ergodicity
is formulated as follows. Let $(\ga,\a)$ be a $C^{*}$--dynamical
system consisting of the unital $C^{*}$--algebra $\ga$ and the automorphism
$\a$. Then the unique ergodicity for $(\ga,\a)$ is equivalent (cf.
\cite{AD}) to the norm convergence
\begin{equation}
\label{er}
\lim_{n\to+\infty}\frac1n\sum_{k=0}^{n-1}\a^{n}(a)=E(a)\,,\quad
a\in\ga\,.
\end{equation}

Here, $E=\om(\,\cdot\,)\idd$ is the conditional expectation onto the fixed--point subalgebra of $\a$
consisting of the constant multiples of the identity, and
$\om\in\cs(\ga)$ is the unique invariant state for $\a$. 
A natural generalization requires that the the fixed--point subalgebra for $\a$ in \eqref{er}
is nontrivial. 
This property, denoted as
the unique ergodicity with respect to the fixed--point subalgebra, has
been investigated in \cite{AD}.
The strict weak
mixing was investigated in \cite{FM}. This means that
\begin{equation*}
\lim_{n\to+\infty}\frac{1}{n}\sum_{k=0}^{n-1}\big|\f(\a^{k}(a))-\f(E(a))\big|=0\,,
\quad a\in\ga\,,
\end{equation*}
for every $\f\in\cs(\ga)$. As before, $E$ is the unique conditional expectation
projecting onto the fixed--point subalgebra. Notice that the last mentioned ergodic property is implied (cf. \cite{Z}) by the following norm convergence 
\begin{equation}
\label{cmcz}
\lim_{n\to+\infty}\frac1n\sum_{k=0}^{n-1}\a^{n_k}(a)=E(a)\,,\quad
a\in\ga
\end{equation}
provided $\{n_k\}_{k\in\bn}$ is any subsequence of natural numbers with nonnull lower density.

In the present paper we investigate another very strong ergodic property as follows. 
We simply require that
\begin{equation}
\label{mta2}
\lim_{n\to+\infty}\f(\a^{n}(a))=\f(E(a))\,,
\quad a\in\ga\,,
\end{equation}
for every $\f\in\cs(\ga)$.
As for the previous situations concerning ergodicity and weak mixing (cf. Proposition \ref{umr}), the convergence to the equilibrium considered here, is connected to the norm convergence of suitable  Cesaro means. Namely, \eqref{mta2} is implied by the norm convergence of the means in \eqref{cmcz},
for all the subsequences $\{n_k\}_{k\in\bn}$ of natural numbers. 

The property \eqref{mta2} of the convergence to the equilibrium is perfectly meaningful in the quantum setting but it has no counterpart in the classical case. Here, there are the main results of the present paper. If a classical system $(X,T)$ fulfils \eqref{mta2} with $E(f)=\int f\di\m$, the support of the unique invariant measure $\m$ is a singleton, that is, it is conjugate to the trivial one--point dynamical system. The last result holds true under the additional condition of separability, that is when $X$ is a compact metric space. On the other hand, we can exhibit some interesting examples of $E$--mixing $C^*$--dynamical system by passing to  the quantum case. Indeed, we show that the shifts on the reduced 
$C^{*}$--algebras of RD--groups, including the free shift on 
the free 
group on infinitely many generators, and the amalgamated free product
$C^{*}$--algebras, enjoy the property of the convergence to the equilibrium \eqref{mta2}. Such a result was recently extended in \cite{DF} to the case of $q$--commutation relations. Namely,  the shift on the $C^*$--algebras generated by the Fock representation of the $q$--commutation
relations has the strong ergodic property \eqref{mta2} (denoted as unique mixing when there is only one invariant state for the dynamics), when $|q|<1$. Thus, we provide nontrivial examples of uniquely mixing $C^*$--dynamical systems for which  the unique invariant state is 
not faithful (case of the $C^{*}$--algebra generated by the $q$--commutation relations), or when it is faithful (case of the $C^{*}$--algebra generated by the self--adjoint part of the generators of the 
$q$--commutation relations for which the restriction of the Fock vacuum is a faithful trace). However, for all these cases, the associated GNS covariant representation is faithful. 

\section{terminology, notation and basic results}

In the present paper we always deal with unital $C^*$ algebras $\ga$ with the identity $\idd$. 
Let  $\ga$ be a $C^*$--algebra. In view of applications to physical systems with nontrivial superselection structure, it is natural to consider different folia $\cf\subset\ga^*$. Fix a representation $\pi$ on the
Hilbert space $\ch_\pi$, and consider the linear space $L^1(\ch_\pi)$
made of all the trace class operators acting on $\ch_\pi$.
Define for $A\in L^1(\ch_\pi)$, $\f_A(x):=\tr(A\pi(x))$, and 
$\cf_\pi:=\{\f_A\,|\,A\in L^1(\ch_\pi)\}$.
Let $\cf\subset\ga^*$ be a set of linear functionals of a $C^*$--algebra $\ga$.
It is called a {\it folium} if it arises as $\cf=\cf_\pi$, for some
representation $\pi$ of $\ga$.

For a (discrete) $C^*$--dynamical system we mean a triplet $\big(\ga,\a,\om\big)$ 
consisting of a unital $C^*$-algebra $\ga$, an automorphism $\a$ of $\ga$, 
and a state $\om\in\cs(\ga)$ invariant under the action of $\a$. The pair $(\ga,\a)$ consisting of 
$C^*$-algebra and an automorphism $\a$ as before, is called a $C^*$--dynamical system as well. A classical
$C^*$--dynamical system is simply a dynamical system such that $\ga\sim C(X)$, $C(X)$ being the Abelian $C^*$-algebra of all the continuous functions on the compact Hausdorff space $X$. In this situation, $\a(f)=f\circ T$ for some homeomorphism $T:X\mapsto X$. 

A $C^*$--dynamical system $\big(\ga,\a\big)$ is said to be {\it uniquely ergodic} if there exists only one state $\om$ invariant for $\a$. A classical $C^*$--dynamical system $(C(X), T)$ is said to be {\it strictly ergodic}, if it is uniquely ergodic and $\m(U)>0$ for each open set $U\subset X$, $\m$ being the unique invariant probability measure invariant under $T$. $(C(X), T)$ is said to be {\it minimal} if 
$\overline{\{T^nx\,|\,\n\in\bz\}}=X$ for each $x\in X$. It is said {\it forward topological transitive} if 
$\overline{\{T^nx\,|\,\n\in\bn\}}=X$ for some $x\in X$.
Let $(X_j,\ca_j,\m_j,T_j)$ be measurable dynamical systems consisting for $j=1,2$, of sets $X_j$, 
$\s$--algebras 
$\ca_j$, probability measures $\m_j$ on $\ca_j$,
and finally measure--preserving invertible transformations $T_j$. The last are said to be {\it conjugate} if there exist sets $A_j\in\ca_j$ of full measure such that 
$T_j(A_j)=A_j$, and a one--to--one measure--preserving map $S:A_1\mapsto A_2$ such that
$T_2=S\circ T_1\circ S^{-1}$.

Let $\big(\ga,\a,\om\big)$ be a $C^*$--dynamical system  and $a,b\in\ga$.
It is said to be {\it ergodic}, {\it weakly mixing} or {\it mixing} if
\begin{align*}
&\lim_{n\to+\infty}\frac{1}{n}\sum_{k=0}^{n-1}\om(a\a^{k}(b))=
\om(a)\om(b)\,,\\
&\lim_{n\to+\infty}\frac{1}{n}\sum_{k=0}^{n-1}\big|\om(a\a^{k}(b))-
\om(a)\om(b)\big|=0\,,\\
&\lim_{n\to+\infty}\om(a\a^{n}(b))=\om(a)\om(b)\,,
\end{align*}
respectively.
Let $\big(\ch_{\om},\pi_{\om},U,\Om\big)$ be the GNS covariant representation 
canonically associated to the dynamical system under 
consideration. It can be straightforwardly seen (see e.g. Proposition \ref{ppp} for a similar situation)
that $\big(\ga,\a,\om\big)$ is ergodic, respectively 
weakly mixing or mixing if and only if
\begin{equation}
\label{mxmx1}
\lim_{n\to+\infty}\frac{1}{n}\sum_{k=0}^{n-1}\langle U^{k}\xi,\eta\rangle=
\langle\xi,\Om\rangle\langle\Om,\eta\rangle\,,
\end{equation}
\begin{equation}
\label{mxmx2}
\lim_{n\to+\infty}\frac{1}{n}\sum_{k=0}^{n-1}\big|\langle U^{k}\xi,\eta\rangle
-\langle\xi,\Om\rangle\langle\Om,\eta\rangle\big|=0\,,
\end{equation}
\begin{equation}
\label{mxmx3}
\lim_{n\to+\infty}\langle U^{n}\xi,\eta\rangle=
\langle\xi,\Om\rangle\langle\Om,\eta\rangle\,,
\end{equation}
respectively.
Let $s(\om)$ be the support of $\om$ in the bidual 
$\ga^{**}$. Then $s(\om)\in Z(\ga^{**})$ if and only if $\Om$ is 
also separating for $\pi(\ga)''$, $Z(\ga^{**})$ being the centre of 
$\ga^{**}$ (see e.g. \cite{SZ}, Section 10.17).
Recall that a {\it conditional expectation}
$E:\ga\mapsto\gb\subset\ga$ is a norm--one projection of the
$C^{*}$--algebra $\ga$ onto a $C^{*}$--subalgebra $\gb$.  

Let $(\ga,\a)$ be a $C^*$--dynamical system, and $E:\ga\mapsto\ga$ a linear map.
\begin{defin}${}$\\
\label{smx}
\begin{itemize}
\item[(i)] $(\ga,\a)$ is said to be
$E$--ergodic if
\begin{equation*}
\lim_{n\to+\infty}\frac{1}{n}\sum_{k=0}^{n-1}\f(\a^{k}(x))=\f(E(x))\,,\quad x\in\ga\,,\f\in\cs(\ga)\,.
\end{equation*}
\item[(ii)] $(\ga,\a)$ is said to be
$E$--weakly mixing if
\begin{equation*}
\lim_{n\to+\infty}\frac{1}{n}\sum_{k=0}^{n-1}\big|\f(\a^{k}(x))-\f(E(x))\big|=0\,,\quad
x\in\ga\,,\f\in\cs(\ga)\,.
\end{equation*}
\item[(iii)] $(\ga,\a)$ is said to be
$E$--mixing if
\begin{equation}
\label{mppp1}
\lim_{n\to+\infty}\f(\a^{n}(x))=\f(E(x))\,,\quad x\in\ga\,,\f\in\cs(\ga)\,.
\end{equation}
\end{itemize}
\end{defin}

It can readily seen (cf. \cite{FM}) that the map $E$ is a conditional expectation 
projecting onto the fixed point subalgebra $\ga^{\a}$.
If $E=\om(\,\cdot\,)\idd$ (i.e. when there is a unique invariant state for $\a$), we call the dynamical sistem 
under consideration {\it uniquely
ergodic}, {\it uniquely weak mixing} or {\it uniquely mixing} respectively.
By using the Jordan decomposition of bounded linear
functionals, one can replace
$\cs(\ga)$ with $\ga^*$ everywhere in Definition \ref{smx}.
\begin{prop}
\label{uesm} 
Let $\big(\ga,T\big)$ be a $C^*$-dynamical system. Then
(iii)$\Rightarrow$(ii)$\Rightarrow$(i).
\end{prop}
\begin{pf}
Let $\f\in\cs(\ga)$. 
\begin{align*}
\text{(iii)}\Rightarrow\text{(ii)}\,\,&\f(\a^{k}(x))\to\f(E(x))\iff |\f(\a^{k}(x))-\f(E(x))|\to0\\
&\Rightarrow\frac{1}{n}\sum_{k=0}^{n-1}\big|\f(\a^{k}(x))-\f(E(x))\big|\to0\,.\\
\text{(ii)}\Rightarrow\text{(i)}\,\,&
\bigg|\frac{1}{n}\sum_{k=0}^{n-1}\big(\f(\a^{k}(x))-\f(E(x))\big)\bigg|
\leq\frac{1}{n}\sum_{k=0}^{n-1}\big|\f(\a^{k}(x))-\f(E(x))\big|\to0
\end{align*}
if (ii) holds true.
\end{pf}
 
According to the case of $E$--ergodicity (cf. \cite{AD}) and $E$--weak mixing (cf. \cite{Z}), we provide a sufficient condition for the $E$--mixing concerning the norm convergece of Cesaro means  of all the subsequences $\{\a^{n_k}\}_{k\in\bn}$.
\begin{prop}
\label{umr}
Let $(\ga,\a)$ be a $C^*$--dynamical system. Suppose that there exists a conditional expectation
$E:\ga\mapsto\ga^{\a}$ of $\ga$ onto the fixed--point subalgebra $\ga^{\a}$,
such that 
\begin{equation}
\label{bks}
\lim_{n\to+\infty}\frac{1}{n}\bigg\|\sum_{k=1}^{n}\a^{n_k}(x)\bigg\|=0
\end{equation}
for each $x\in\ga$ satisfying $E(x)=0$, and for each sequence $0\leq n_1<n_2<\cdots<n_k<\cdots$ of increasing natural numbers. 
Then $(\ga,\a)$ is $E$--mixing.
\end{prop}
\begin{pf}
Suppose that \eqref{bks} holds true for each sequence 
$0\leq n_1<n_2<\cdots<n_k<\cdots$ of increasing natural numbers and for each $x\in\ga$ satisfying 
$E(x)=0$, but $(\ga,\a)$ is not $E$--mixing. Then
there should exists a norm one
functional $\f\in\ga^*$ and an element $x\in\ga$ with $E$--vanishing expectation, such that
$\f(\a^n(x))$ does not vanish when $n\to+\infty$. By passing to subsequences and
modifying $\f$ by a phase--factor, we can suppose that $\re(\f(\a^{n_k}(x))\geq c>0$,
$k=0,1,2,\dots$\,\,. 
But
\begin{equation*}
\bigg\|\frac1n\sum_{k=0}^{n-1}\alpha^{n_k}(x)\bigg\|
\geq\frac1n\sum_{k=0}^{n-1}\re(\f(\a^{n_k}(x))\geq c>0
\end{equation*}
which contradicts \eqref{bks}. We then conclude that ${\displaystyle\lim_{n\to+\infty}\f(\a^{n}(x))=0}$
for each $x$ with $E$--vanishing expectation. The result follows as $a\in\ga$ can be written as
$a=(a-E(a))+E(a)$ with $E(a)\in\ga^{\a}$. 
\end{pf}

\section{strong mixing of the free shifts}
\label{qddss}

The present section is devoted to prove that some quantum dynamical systems enjoy
the very strong ergodic property of $E$--mixing, $E$ being the
conditional expectation on the fixed point subalgebra. The examples under considerations 
are
the free shift on the reduced
amalgamated free product $C^{*}$--algebra, and  lenght--preserving
automorphisms of
the reduced $C^*$--algebra of RD--group for the lenght--function, the latter 
including the free shift on the free group on infinitely many
generators. 

Let $\gd$ be a unital $C^{*}$--algebra with identity $\idd$, and
$E^{\gd}_{\gb}:\gd\mapsto\gb$ a conditional expectation onto the unital
$C^{*}$--subalgebra $\gb$ with the same identity $\idd$. For each
integer $i\in\bz$, consider a copy $(\ga_i,E_i)$ of $(\gd,E^{\gd}_{\gb})$,
together with the reduced amalgamated free product
\begin{equation}
\label{eq:Aphi} 
(\ga,E)={(*_\gb)}_{i\in\bz}(\ga_i,E_i)\,.
\end{equation}

The $C^{*}$--algebra $\ga$ naturally acts on a Hilbert right
$\gb$--module $\ce$ and it is generated by 
$\big\{\l^{i}_{a}\,:\, a\in\ga_{i}\,,i\in\bz\big\}$, $\l^{i}$ being the embedding of $\ga_{i}$ in
$\cb_{\gb}(\ce)$, the space of all the bounded $\gb$--linear maps acting
on $\ce$.
The conditional expectation $E$ is given by
$$
E(a)=\langle\idd a\,,\,\idd\rangle\,,\qquad a\in\ga\,,
$$
$\langle\,\cdot\,,\,\cdot\,\rangle$ being the $\gb$--valued inner
product of $\ce$ which is supposed be linear w.r.t. the first
variable. We refer the reader to \cite{AD} and the references cited therein, for further details.

The free--shift automorphism $\alpha$ on $\ga$ is the automorphism
of $\ga$ given by $\alpha(\lambda_a^i)=\lambda_a^{i+1}$ for all
$a\in \ga$ and $i\in\bz$.
\begin{thm}
\label{F-mix}
Let $\alpha$ be the free--shift automorphism on the
reduced amalgamated free product $C^*$--algebra $\ga$ given in
\eqref{eq:Aphi}. Then $\alpha$ is $E$--mixing.
\end{thm}
\begin{pf}
Fix $a\in \ga$ of the form $a=w$ for a word
$w=\l^{m(1)}_{a_1}\l^{m(2)}_{a_2}\cdots\l^{m(p)}_{a_p}$, with
$p\geq 1$, $a_i\in \ga_{m(i)}^\circ$, and $m(i)\in\bz$ fulfilling
$m(i)\ne m(i+1)$, $i=1,\dots,p-1$. Here if $i\in\bz$,
$\ga_{i}^\circ:=\big\{a-E_{i}(a)\,:\,a\in \ga_{i}\big\}$.
Notice that
$$
\a^{k}(w)=\l^{m(1)+k}_{a_1}\l^{m(2)+k}_{a_2}\cdots\l^{m(p)+k}_{a_p}\,.
$$
Let us take any increasing sequence $\{k_j\}_{j\in\bn}\subset\bn$.
Define
$$
m_{j}(l):=m(l)+k_{j}\,,\quad j\in\bn\,,\,\,\,1\leq l\leq p\,.
$$
Notice that
$$
\a^{k_{j}}(w)=\l^{m_{j}(1)}_{a_1}\l^{m_{j}(2)}_{a_2}
\cdots\l^{m_{j}(p)}_{a_p}\,.
$$
In addition, $m_{j}(1)\neq m_{j}(p-1)$, $j\in\bn$, and $j\neq j'$
implies $m_{j}(1)\neq m_{j'}(1)$, $m_{j}(p)\neq m_{j'}(p)$.
Thus, we can apply the estimation in Proposition 5.1 of \cite{AD}
to the element ${\displaystyle f_{n}:=\sum_{j=1}^{n}\a^{k_{j}}(w)}$,
obtaining
\begin{equation}
\label{hinq}
\bigg\|\frac{1}{n}\sum_{j=1}^{n}\alpha^{k_j}(w)\bigg\|
\leq\frac{2p+1}{n^{1/2}}\prod_{i=1}^p\|a_i\|\,.
\end{equation}
It was proven in \cite{AD} that $\a$ is $E$--uniquely ergodic.
This implies that
\begin{equation*}
\lim_{n\to+\infty}\frac1n\sum_{k=0}^{n-1}\alpha^k(a)=E(a)
\end{equation*}
for every $a\in\ga$.
By a standard density argument, \eqref{hinq} implies that, if $a\in\ga$ fulfils $E(a)=0$, then
$$
\lim_{n\to+\infty}\bigg\|\frac{1}{n}\sum_{j=1}^{n}\alpha^{k_j}(a)\bigg\|=0\,.
$$
The proof is now complete by taking into account Proposition \ref{umr}.
\end{pf}

The other examples considered here are the free shifts on the
reduced $C^*$--algebra on RD--groups (cf.  \cite{J}).
\begin{thm}
\label{ue1} 
Let $\b$ be a lenght--preserving automorphism of a
RD--group $G$ for the lenght--function $L$, such that its orbits
are infinite or singletons. Then the automorphism $\a$ induced by
$\b$ on $C^{*}_{r}(G)$ is $E$--mixing.
\end{thm}
\begin{pf}
Let $H:=\{g\in G\,:\,\b(g)=g\}$. As $\a$ is $E$--uniquely ergodic
(cf. \cite{AD}, Proposition 3.5), the pointwise limit in norm
$$
E:=\lim_{n\to+\infty}\frac{1}{n}\sum_{k=1}^{n}\a^{k}
$$
exists and gives rise to the conditional expectation projecting onto
the fixed--point algebra $C^{*}_{r}(H)\subset C^{*}_{r}(G)$. The
proof follows as in
Theorem \ref{F-mix} by taking into account that
$$
\lim_{n\to+\infty}\bigg\|\frac{1}{n}\sum_{j=1}^{n}\a^{k_{j}}(\l_{g})\bigg\|
\leq C(1+L(g))^{s}\bigg\|\frac{1}{n}\sum_{j=1}^{n}\d_{\b^{k_{j}}(g)}
\bigg\|_{\ell^{2}(G)}=\frac{C(1+L(g))^{s}}{\sqrt{n}}
$$
for each sequence $\{k_{j}\}$ of natural numbers, and $\b(g)\neq g$.
\end{pf}

Finally, we report the case of the automorphism generated by the
shift on the free group on infinitely many generators. The shift on the generators is defined as
$\b:g_i\mapsto g_{i+1}$, $i\in\bz$.
\begin{cor}
Let $\bbf_{\infty}$ be the free group on infinitely many
generators $\{g_{i}\}_{i\in\bz}$. Then the automorphism $\a$ induced on
$C^{*}_{r}(\bbf_{\infty})$ by the free shift on the generators is
$E$--mixing with $E=\t(\,\cdot\,)\idd$, $\t$ being
the canonical trace on $C^{*}_{r}(\bbf_{\infty})$.
\end{cor}
\begin{pf}
By taking into account the Haagerup inequality (cf. \cite{H},
Lemma 1.4), we reduce the matter to a particular case
of both Theorem \ref{F-mix} and Theorem \ref{ue1}.
\end{pf}

\section{the case of classical dynamical systems}
\label{epcs}

Let $(\ga,\a)$ be a $C^*$--dynamical system. Consider $E$--ergodicity, $E$--weak mixing and $E$--mixing listed before. When there exists only one invariant state $\om$ for $\a$ (i.e. when 
$E=\om(\,\cdot\,)\idd$), the last can be considered as the topological
analogous of measure theoretic ergodicity, weakly mixing and mixing, respectively. 
It is known that the irrational rotations on the circle satisfy
(i) but not (ii) in Definition \ref{smx}. However, at least for the separable case, we will show that it is impossible to exhibit a classical dynamical system satisfying (iii) in Definition \ref{smx}, unless it is conjugate to the trivial one--point dynamical system. 

In the present section we restrict our analysis to compact metric spaces.
To have an idea of what happens in the classical situation, we discuss the following example  
suggested by D. Kerr. We start with
the one--point compactification $\bz_\infty$ of $\bz$. The
shift $\a(f)(x)=f(x+1)$ extends to an automorphism of $C(\bz_\infty)$,
and the measure defined by the evaluation at infinity
$$
\m_\infty(f):=\lim f(x)\equiv f(\infty)
$$
is invariant for $\a$. 
The dynamical system under consideration satisfies 
\eqref{mppp1}, but it is not strictly ergodic. However, it is conjugate to the trivial one--point dynamical system. We prove that all the classical dynamical system satisfying \eqref{mppp1}, arise essentially in this way, up to conjugacy.\footnote{Other similar examples come from the one point compactification of countably many copies of $\bz$, by following the previous construction.} We start with the following preparatory result.
\begin{lem}
\label{prre}
Let $(X, T)$ be strictly ergodic. Then there exists some $x_0\in X$ and a subsequence 
$\{n_x(k)\}_{k\in\bn}$ such that ${\displaystyle\lim_k T^{n_x(k)}x_0=x}$ for each $x\in X$.
\end{lem}
\begin{pf}
Let $d$ be a fixed metric on $X$. We can restrict the matter to the case $|X|=+\infty$, the result being trivial for finite $X$. Then $X$ is a perfect set, otherwise it would be finite. The dynamical system
$(X,T)$ is minimal (cf. \cite{W}). Thus, it is forward topological transitive, 
see e.g. \cite{Wa}, Theorem 5.10 by taking account the remark after Theorem 5.6. Then there exists a forward dense orbit $\{T^nx_0\,|\,n\in\bn\}$ for some $x_0\in X$. Fix $x\in X$. We can find a function 
$f_x:\bn\to\bn$ such that
\begin{equation}
\label{iiff}
0<d(x,T^{f_x(l)}x_0)<1/l\,,\qquad l=1,2,\dots\,\,.
\end{equation}

Notice that the range of $f_x$ cannot be finite. Indeed, if $|f_x|<+\infty$, there would exist infinitely many $l$, say $\{l_k\}_{k\in\bn}$, such that for some integer $n_0$,
$$
f_x(l_k)=n_0\,,\qquad l=1,2,\dots\,\,.
$$
Thus, by \eqref{iiff},
$$
0<d(x,T^{n_0}x_0)<1/l_k\to0\,.
$$

As the sequence $\{f_x(l)\}_{l\in\bn}$ of natural numbers has $+\infty$ as a cluster point, there exists a subsequence, defined as $n_x(k):=f_x(l_k)$, such that ${\displaystyle\lim_k n_x(k)}=+\infty$. Finally, by construction, ${\displaystyle\lim_k T^{n_x(k)}x_0=x}$.
\end{pf}
\begin{thm}
\label{z3}
Let $X$ be a compact metric space, and $T:X\mapsto X$ a homeomorphism.
Suppose that $(X, T)$ satisfies \eqref{mppp1}, with $\m$ the unique invariant probability measure under $T$ and $E(f)=\int f\di\m$.

Then $(X, T)$ is conjugate to the trivial one--point dynamical system.
\end{thm}
\begin{pf}
Consider the  support $\supp(\m)$ of the unique invariant measure for the trasformation $T$. It is well--known that $\supp(\m)$ is a closed
invariant set for $T$. Then the resulting dynamical system
$\big(\supp(\m),T\lceil_{\supp(\m)}\big)$ is strictly ergodic by construction. By Lemma \ref{prre},
there exists $x_0\in\supp(\m)$, and a subsequence $\{n_x(k)\}_{k\in\bn}$ depending on $x\in\supp(\m)$,
such that ${\displaystyle\lim_kT^{n_{x}(k)}x_0=x}$. Fix
$x_1,x_2\in\supp(\m)$, and consider $f\in C(\supp(\m))$. Then thanks to \eqref{mppp1},
$$
f(x_1)=\lim_kf(T^{n_{x_1}(k)}x_0)=\m(f)=\lim_kf(T^{n_{x_2}(k)}x_0)=f(x_2)\,.
$$
This means that $x_1=x_2$.
\end {pf}

According to the Jewett--Krieger theorem (cf. \cite{J, K}), it is possible to construct a huge class of uniquely ergodic classical dynamical systems. It is still unclear if the Jewett--Krieger theorem can be established in the weak mixing situation. However, by Theorem \ref{z3}, the analogous of the Jewett--Krieger theorem cannot be carried out for the mixing situation.
\begin{cor}
\label{zz4}
Let $(X,\ca,\m,T)$ be a measurable dynamical system. If it is conjugate to a uniquely mixing one, then $L^{\infty}(X,\ca,\m)\sim\bc\idd$.
\end{cor}
\begin{pf}
Theorem  \ref{z3} implies that if $(X,\ca,\m,T)$ is conjugate to a uniquely mixing one, then it is conjugate to the trivial one--point dynamical system as well.
\end{pf}

\section{more on dynamical systems}

 At the light of the previous results, it is natural to address some ergodic properties which are suitable for applications to quantum physics.

Let $(\ga,\a)$ be a $C^*$--dynamical system, and $\om\in\cs(\ga)$. Consider for  $a,b,c\in\ga$, the following properties
\begin{equation}
\label{mmmm1}
\lim_{n\to+\infty}\frac1n\sum_{k=0}^{n-1}\om(b\a^k(a)c)=\om(bc)\om(a)\,,
\end{equation}
\begin{equation}
\label{mmmm2}
\lim_{n\to+\infty}\frac1n\sum_{k=0}^{n-1}\big|\om(b\a^k(a)c)-\om(bc)\om(a)\big|=0\,,
\end{equation}
\begin{equation}
\label{mmmm3}
\lim_{n\to+\infty}\om(b\a^n(a)c)=\om(bc)\om(a)\,.
\end{equation}
Notice that by polarization,
\eqref{mmmm1},\eqref{mmmm2}, \eqref{mmmm3} are equivalent to
\begin{align*}
&\lim_{n\to+\infty}\frac1n\sum_{k=0}^{n-1}\om(b^*\a^k(a)b)=\om(b^*b)\om(a)\,,\\
&\lim_{n\to+\infty}\frac1n\sum_{k=0}^{n-1}\big|\om(b^*\a^k(a)b)-\om(b^*b)\om(a)\big|=0\\
&\lim_{n\to+\infty}\om(b^*\a^n(a)b)=\om(b^*b)\om(a)\,,
\end{align*}
for each $a,b\in\ga$, respectively.

Properties \eqref{mmmm1}, \eqref{mmmm2}, \eqref{mmmm3} imply 
ergodicity, weak mixing or mixing for the dynamical system $(\ga,\a,\om)$, respectively. It is expected that the former are stronger then the latter.
However we discuss relevant situations for which we have the equivalence between the corresponding properties. 

We start by saying that the
$C^*$ dynamical system $(\ga,\a,\om)$ is 
{\it asymptotically Abelian} if
$$
\lim_{n\to\pm\infty}\om(c[\a^n(a),b]d)=0\,,\quad a,b,c,d\in\ga\,,
$$
where $[a,b]\equiv ab-ba$ stands for the commutator.\footnote{The reader is referred to \cite{F} and the literature cited therein, for the more general case when Fermions are present. Proposition \ref{ewmm} holds true as well if the invariant state$\om$ is graded asymptotically Abelian.}

The following result is more or less known to the experts. We report it for the sake of completeness.
\begin{prop}
\label{ewmm}
Let the $C^*$ dynamical system $(\ga,\a,\om)$ be 
asymptotically Abelian, or the support of $\om$ be central in $\ga^{**}$. Then  the properties 
\eqref{mmmm1}, \eqref{mmmm2}, \eqref{mmmm3} are equivalent to 
ergodicity, weak mixing or mixing for the dynamical system $(\ga,\a,\om)$, respectively.
\end{prop}
\begin{pf}
Let the $C^*$ dynamical system $(\ga,\a,\om)$ be asymptotically Abelian, and set
$\G_n:=\om(b[\a^n(a),c])$. By asymptotic Abelianess $\G_n\longrightarrow0$, and
$$
\frac1n\sum_{k=0}^{n-1}\G_k\longrightarrow0\,,\qquad
\frac1n\sum_{k=0}^{n-1}|\G_k|\longrightarrow0
$$
as well. Thanks to this, we get for example,
\begin{align*}
&\frac1n\sum_{k=0}^{n-1}\big|\om(b\a^k(a)c)-\om(bc)\om(a)\big|\leq\frac1n\sum_{k=0}^{n-1}|\G_k|\\
+&\frac1n\sum_{k=0}^{n-1}\big|\om(bc\a^k(a))-\om(bc)\om(a)\big|\longrightarrow0
\end{align*}
by weak mixing. The other equivalences follow analogously.

Let now $(\ch_{\om},\pi_{\om},\Om,U)$ be the GNS covariant representation for $(\ga,\a,\om)$.  Suppose that the $C^*$ dynamical system $(\ga,\a,\om)$ is ergodic, weakly mixing or mixing. Then it satisfies \eqref{mxmx1}, \eqref{mxmx2}, \eqref{mxmx3} respectively.
As the remaining assertion follow analogously, we restrict ourselves to mixing situation.
Let $T\in\pi_{\om}(\ga)'$. We have
\begin{equation}
\label{ccy}
\langle\pi_{\om}(b)\pi_{\om}(\a^n(a))T\Om,\Om\rangle\longrightarrow
\langle\pi_{\om}(b)T\Om,\Om\rangle\om(a)\,.
\end{equation}

As the 
support of $\om$ is central in $\ga^{**}$, $\Om$ is cyclic for $\pi_{\om}(\ga)'$ too. Then \eqref{ccy} leads to
\begin{equation}
\label{ccy1}
\langle\pi_{\om}(b)\pi_{\om}(\a^n(a))\xi,\Om\rangle\longrightarrow
\langle\pi_{\om}(b)\xi,\Om\rangle\om(a)
\end{equation}
for generic $\xi\in\ch_{\om}$.
The assertion follows by computing \eqref{ccy1} with $\xi=\pi_{\om}(c)\Om$.
\end{pf}

We start with various convergences to the equilibrium in a fixed folium. These are intermediate situations between ergodicity, weak mixing and mixing (i.e. when one considers the folium generated by a fixed invariant state), and unique ergodicity, unique weak mixing and unique mixing (i.e. if  one considers the universal folium $\ga^*$ in the case when the invariant state is unique).

Let $(\ga,\a)$ be a $C^*$--dynamical system, $\cf\subset\ga^*$ a folium, and $\om\in\cs(\ga)$ a
state.
We say that $(\ga,\a)$ is 
$(\cf,\om)$--{\it ergodic}, $(\cf,\om)$--{\it weakly mixing} or
$(\cf,\om)$--{\it mixing} if
\begin{equation*}
\lim_{n\to+\infty}\frac{1}{n}\sum_{k=0}^{n-1}\f(\a^{k}(x))=\f(\idd)\om(x)\,,
\end{equation*}
\begin{equation*}
\lim_{n\to+\infty}\frac{1}{n}\sum_{k=0}^{n-1}\big|\f(\a^{k}(x))-\f(\idd)\om(x)\big|=0\,,
\end{equation*}
\begin{equation*}
\lim_{n\to+\infty}\f(\a^{n}(x))=\f(\idd)\om(x)\,.
\end{equation*}
holds true for each $x\in\ga$, $\f\in\cf$ respectively.
\begin{prop}
\label{ppp}
Let  $(\ga,\a,\om)$ be a $C^*$--dynamical system. Then
$(\cf_{\pi_{\om}},\om)$--ergodicity, $(\cf_{\pi_{\om}},\om)$--weak mixing or 
$(\cf_{\pi_{\om}},\om)$--mixing are equivalent to
\eqref{mmmm1}, \eqref{mmmm2}, \eqref{mmmm3}, respectively. 
\end{prop}
\begin{pf}
We have only to show that each of the latter properties implies the corresponding one in the former list. We restrict ourselves to weak mixing, as the remaining ones follow analogously. 
Let $\f\in\cf_{\pi_{\om}}$ and $x\in\ga$. Then there exists 
${\displaystyle\xi,\eta\in\bigoplus_{\bn}\ch_{\om}}$ such that
${\displaystyle\f(x)=\langle\bigoplus_{\bn}\pi_{\om}(x)\xi,\eta\rangle}$.
If $\e>0$, then there exist $L\in\bn$ such that
${\displaystyle\bigg|\sum_{l>L}\langle\pi_{\om}(x)\xi_l,\eta_l\rangle\bigg|<\e\|x\|}$. 
Put ${\displaystyle K:=1+\|\eta\|+\sqrt{1+2\|\xi\|+\|\xi\|^2}}$.
Choose $b_l$, $c_l$, such that
${\displaystyle\|\xi_l-\pi_{\om}(c_l)\Om\|<\e/\sqrt{L}\,,\quad\|\eta_l-\pi_{\om}(b_l^*)\Om\|<\e/\sqrt{L}}$,
$l=1,\dots,L$, and finally fix $a\in\ga$ with $\|a\|\leq1$. 
It is straightforward to check that
$$
\frac{1}{n}\sum_{k=0}^{n-1}\big|\f(\a^{k}(a))-\f(\idd)\om(a)\big|\leq2K\e
+\sum_{l=1}^{L}
\bigg(\frac{1}{n}\sum_{k=0}^{n-1}\big|\om(b_l\a^k(a)c_l)-\om(b_lc_l)\om(a)\big|\bigg)
$$
which goes to zero as $\e$ is arbitrary.
\end{pf}

Notice that, the $(\ga^*,\om)$--ergodicity, 
$(\ga^*,\om)$--weak mixing, or $(\ga^*,\om)$--mixing mean by definition, unique ergodicity, unique weak mixing and unique mixing respectively, $\om$ being the unique invariant state. 
The $(\cf_{\pi_{\om}},\om)$--ergodicity, 
$(\cf_{\pi_{\om}},\om)$--weak
mixing or $(\cf_{\pi_{\om}},\om)$--mixing are nothing but, the natural generalizations suitable for physical applications (cf. Proposition 
\ref{ewmm}), of the standard ergodicity, weak mixing and mixing respectively, 
for the dynamical system $(\ga,\a,\om)$. 
It is easy to show that
$\om$ is invariant, but not necessarily
$\om\in\cf$. In addition, it is unknown if the forward ergodic properties as those listed
above, do imply the corresponding ones for the backward dynamics.
However, we have
\begin{prop}
Let $(\ga,\a)$ be a $C^{*}$--dynamical system.
\begin{itemize}
\item[(i)]
If $(\ga,\a)$ is $(\ga^*,\om)$--ergodic, $(\ga^*,\om)$--weakly
mixing, or $(\ga^*,\om)$--mixing, then
$(\ga,\a^{-1})$ is $(\ga^*,\om)$--ergodic. 
\item[(ii)] If $(\ga,\a)$ is $(\cf_{\pi_{\om}},\om)$--ergodic,
$(\cf_{\pi_{\om}},\om)$--weakly
mixing or $(\cf_{\pi_{\om}},\om)$--mixing, then
$(\ga,\a^{-1})$ enjoys the corresponding property, provided $(\ga,\a,\om)$ satisfies one of the hypotheses of Proposition \ref{ewmm}.
\end{itemize}
\end{prop}
\begin{pf}
(i) If $(\ga,\a)$ fulfils anyone of the properties listed above, it is
uniquely ergodic, with $\om$ as the unique invariant state. But $\om$
is the unique invariant state for $\a^{-1}$ as well. This means that
$(\ga,\a^{-1})$ is $(\ga^*,\om)$--ergodic.

(ii) Let $U$ be the unitary implementing $\a$ on $\ch_{\pi_{\om}}$.
According to Proposition \ref{ewmm}, the dynamical system $(\ga,\a,\om)$ is ergodic, weakly mixing or mixing if and
only if \eqref{mxmx1}, \eqref{mxmx2} or \eqref{mxmx3} holds true respectively. It is easy to show that each one of the latter is satisfied if we replace $U$ with $U^*$. The assertion will follow again by Proposition \ref{ewmm}, as
the canonical implementation of $\a^{-1}$ is precisely $U^{*}$.  
\end{pf}

It would be interesting to
construct a $C^*$--dynamical system $(\ga,\a)$ which is
uniquely weakly mixing or uniquely mixing, such that
$(\ga,\a^{-1})$ is not, provided the such a dynamical system would exist. This is not the case of the free shift on the free group.
\begin{prop}
The dynamical system $(C^{*}_{r}(\bbf_{\infty}),\a^{-1})$ is 
$(C^{*}_{r}(\bbf_{\infty})^*,\t)$--mixing as well. 
\end{prop}
\begin{pf}
Let the $g_{k}$ be the $k$--generator of 
$\bbf_{\infty}$, and $\th\l_{g_k}:=\l_{g_{-k}}$ he "time reversal" symmetry.\footnote{Recall that a time reversal symmetry does not always exist for a dynamical system, see e.g \cite{SW}.}
Such an automorphism $\th$ fulfils $\th^2=\id$, $\th\a=\a^{-1}\th$, $\t\circ\th=\t$. We get for $x\in C^{*}_{r}(\bbf_{\infty})$, $\f\in C^{*}_{r}(\bbf_{\infty})^*$,
\begin{align*}
&\f(\a^{-n}(x))=(\th^*\f)(\th\a^{-n}(x))=(\th^*\f)(\a^n(\th x))\\
&\longrightarrow(\th^*\f)(\idd)(\th^*\t)(x)=\f(\idd)\t(x)\,.
\end{align*}
\end{pf}

\section*{acknowledgment} 
\noindent
We thank K. Dykema and D. Kerr for some useful discussions.


\begin{thebibliography}{9999}



\bibitem{AD} Abadie B., Dykema K.
{\it Unique ergodicity of free shifts and some other automorphisms
of $C^*$--algebras}, J. Operator Theory {\bf 61} (2009), 279--294.

\bibitem{DF} Dykema K., Fidaleo F. {\it Unique mixing of the shift on the $C^*$--algebras generated by the $q$--canonical commutation relations}, Houston J. Math., to appear.

\bibitem{F} Fidaleo F. {\it KMS states and the chemical potential for disordered systems}, Commun. Math. Phys. {\bf 262}  (2006), 373--391.

\bibitem{FM} Fidaleo F., Mukhamedov F. {\it Strict weak mixing of
some $C^*$--dynamical systems based on free shifts},
J. Math. Anal. Appl. {\bf 336} (2007), 180--187.

\bibitem{H} Haagerup U.
{\it An example of a non nuclear $C^{*}$--algebra which has the
metric approximation property}, Invent. Math. {\bf 50} (1979),
279--293.

\bibitem{Je} Jewett R. I. {\it The pervalence of uniquely ergodic
systems}, J. Math. Mec.
{\bf 19} (1970), 717--729.

\bibitem{J} Jolissaint P.
{\it Rapidly decreasing functions in reduced $C^{*}$--algebras of
groups}, Trans. Amer. Math. Soc. {\bf 317} (1979), 279--293.

\bibitem{K} Krieger {\it On unique ergodicity}, in  "Proceedings of
the Sixth Berkeley Symposium on Mathematical Statistics and
Probability (1970/1971), Vol. II: Probability theory". Univ.
California Press, Berkeley, 1972, 327--346.

\bibitem{SZ} Str\v{a}til\v{a} S., Zsid\'o L.
{\it Lectures on von Neumann algebras}, Abacus press, Tunbridge 
Wells, Kent 1979.

\bibitem{SW} Streater R. F., Wightman A. S.
{\it PCT, spin and statistics and all that},   Princeton University Press, New Jersey 2000.

\bibitem{Wa} Walters P.
{\it An introduction to ergodic theory}, Springer,  New York 1982.
Wells, Kent 1979.

\bibitem{W} Weiss B. {\it Strictly ergodic models for dynamical
systems}, Bull. Amer. Math. Soc. (N.S.) {\bf13} (1985), 143--146.

\bibitem{Z}  Zsid\'o L. {\it Weak mixing properties of
vector sequences}, in: The extended field of operator theory,  361--388, Oper. Theory Adv. Appl., 171, BirkhŠuser, Basel, 2007.

\end{thebibliography}
\end{document}